\documentclass[11pt]{article}
   \usepackage{amsfonts}
\usepackage{latexsym,tikz}
\usepackage{amsmath}
 
\title{Multiple Dedekind Zeta Functions}
\author{Ivan Horozov, Nickola Horozov and Zouberou Sayibou}

\newcommand \nc {\newcommand}
\nc \proof {\noindent {\em{Proof.\/ }}} \nc \qed {$\Box$\hfill}
\newtheorem{theorem}{Theorem}[section]
\newtheorem{lemma}[theorem]{Lemma}
\newtheorem{proposition}[theorem]{Proposition}
\newtheorem{corollary}[theorem]{Corollary}
\newtheorem{definition}[theorem]{Definition}
\newtheorem{example}[theorem]{Example}
\newtheorem{remark}[theorem]{Remark}
\newtheorem{conjecture}[theorem]{Conjecture}
\newtheorem{question}[theorem]{Question}
\nc \bth[1] {\begin{theorem}\label{t#1} } \nc \ble[1]
{\begin{lemma}\label{l#1} } \nc \bpr[1]
{\begin{proposition}\label{p#1} } \nc \bco[1]
{\begin{corollary}\label{c#1} } \nc \bde[1]
{\begin{definition}\label{d#1}\rm } \nc \bex[1]
{\begin{example}\label{e#1}\rm } \nc \bre[1]
{\begin{remark}\label{r#1}\rm } \nc \bcon[1]
{\begin{conjecture}\label{con#1}\rm } \nc \bque[1]
{\begin{question}\label{que#1}\rm }
\nc {\eth} { \end{theorem} } \nc {\ele} { \end{lemma} } \nc
{\epr}{ \end{proposition} } \nc {\eco} { \end{corollary} } \nc
{\ede} {\end{definition} } \nc {\eex} { \end{example} } \nc {\ere}
{\end{remark} } \nc {\econ} { \end{conjecture} } \nc {\eque}
{\end{question} }
\def\a{\alpha}

\def\l{\lambda}

\def \Z {{\mathbb Z}}
\def \Q {{\mathbb Q}}
\def \R {{\mathbb R}}
\def \C {{\mathbb C}}

\begin{document}

\title{Distribution of primes represented by polynomials and Multiple Dedekind zeta functions}
\maketitle

\begin{abstract} 
In this paper, we state several conjectures regarding distribution of primes and of  pairs of primes represented by  irreducible homogeneous polynomial in two variables $f(a,b)$. We formulate conjectures with respect to the slope $t=b/a$ for any irreducible polynomial $f$. Here, we formulate a conjecture for all irreducible polynomials.

We also consider conjectures for distribution of pairs of primes. It  show unexpected relation to multiple Dedekind zeta function - 
 at $s=2$ for one prime and at $(s_1,s_2)=(2,2)$ for pairs of primes. 
We tested the conjecture for pairs of primes for several quadratic fields.
The conjecture for pairs of primes and multiple Dedekind zeta function over the Gaussian integers provide error less than a tenth of a percent.

We also tested conjectures that compare sets of primes in a pair of different quadratic fields. Numerically, such quotients can be expressed in terms of regulators and class numbers.

Some of the data, together with the code, is available on GitHub, (see \cite{Zouberou}). 

\end{abstract}



\tableofcontents

\section{Introduction}
Distribution of prime number is a classical problem. There are various ramifications of it, which are quite interesting to this day.
The most striking relation on distribution of primes 
is related to the zeros of the Riemann zeta function. In fact, it is related to the Riemann Conjecture that stating that the  non-trivial zeroes of Riemann zeta function
\[\zeta(s)=\sum_{n=1}^\infty\frac{1}{n^s},\]
where $s$ is a complex number,
are located at the axis $Re(s)=\frac{1}{2}$.
This conjecture is far reaching and has many consequences in number theory and complexity theory.
There are various generalizations of the Riemann zeta function which have corresponding Riemann Hypothesis. 

A central part of the paper consist of conjectures on pairs of primes and their relation to multiple Dedekind zeta values.
We present conjectures on distributions of pairs of primes represented by homogeneous polynomials in two variables.
We use multiple Dedekind zeta functions  \cite{MDZF} to describe distribution of certain pairs of prime numbers.
The conjectures stated here are tested with a very good precision using the programming language Python. 
We have tested them for several real and imaginary quadratic fields.

The simplest version of the conjecture can be stated for distribution of primes (as opposed to a pair of primes) represented by a homogeneous polynomial in two variables and its relations to a multiple Dedekind zeta values.
 We give a conjectural formula for any number field. For real and imaginary quadratic fields this formula is a theorem.

Consider the Gaussian integers $\Z[i]$ consisting of elements $\pi=a+ib$, where $i=\sqrt{-1}$, and $a$ and $b$ are integers. Let $N(\pi)=a^2+b^2$ be the norm of the the Gaussian integer.
We will consider those elements whose norm is a prime number. In other words, we will consider $\pi=a+ib$ such that $N(\pi)=a^2+b^2$ is a prime number. 
Let $a+ib=re^{i\alpha}$ be the polar representation of a complex number.

Let $P(\theta, M)$ be the number of Gaussian integers $\pi=a+ib$ such that 

(1) $N(\pi)=a^2+b^2$ is a prime  number

(2) $N(\pi)<M$

(3) $0<\alpha<\theta$, where $re^{i\alpha}=\pi$.

Due to Hecke \cite{Hecke}, we have the following.
\begin{theorem} 
For large $M$, we have that $P(\theta, M)$  is equidistributed with respect to the angle $\theta$. 
\end{theorem}

This conjecture can be reformulated by expressing $\theta$ in terms of $t=b/a$.
We have $\tan \theta=\frac{b}{a}=t$. Therefore,
\[\theta=\tan^{-1}\left(t\right).\]

For imaginary quadratic fields again, we have equidistribution of primes with respect to the angle $\theta$.
For a negative integer $-D$, we consider the ring of integers $\Z[\sqrt{-D}]$ of the type $a+b\sqrt{-D}$,
where $a$ and $b$ are integers. The norm of $a+b\sqrt{-D}$ is $N(a+b\sqrt{-D})=(a+b\sqrt{-D})(a-b\sqrt{-D})=a^2+Db^2$.

Let $t=b\a$. Then $t\sqrt{D} = \frac{b\sqrt{D}}{a}=\tan\theta$.
Due to 
 Maknys  \cite{M1}, \cite{M2}, \cite{M3}, we have the following.

\begin{theorem} For primes of the form $p =a^2 +Db^2$ 
the distribution for $t=\frac{b}{a}$ is given by inverse tangent
$\tan^{-1}(t\sqrt{D})$. 
\end{theorem}

Then due to Rademacher \cite{R}
we have the following.

\begin{theorem} For primes of the form $p =a^2 -Db^2$ 
 the distribution for $t=\frac{b}{a}$ is given by inverse hyperbolic tangent
 $\tanh^{-1}(t\sqrt{D})$. 
\end{theorem}

\begin{conjecture}
\label{one prime}
If $\alpha$ is a solution of $f(x,1)=0$ where $f(x,y)=x^n+a_1x^{n-1}y+\dots+a_{n-1}xy^{n-1}+a_n$ is is an irreducible polynomial with integer coefficients, then for primes $p$ of the form $p=N(a+b\alpha)$  the  distribution for $t=\frac{b}{a}$ for $p<M$ and $0<t<T$ is $F^{-1}$, where
is 
\[F(T)=c_M\int_0^T\frac{dt}{[\sqrt[n]{f(1,t)}]^2}.\]
\end{conjecture}

{\bf Remarks:}

1. The conjecture is tested for approximately 20 different fields. They are of degree $n$ for $n=3,4,5,6$

2. For cubic field, that is for $n=3$, the integrant is a rational function on the elliptic curve $y^3=f(x)$ where $f$ is an irreducible cubic polynomial.

3. For $n=2$ and $D$ positive, we can solve the integral explicitly.
\[F(T)=\int_T\frac{dt}{1+Dt^2}=\tan^{-1}(t\sqrt{D}).\]

4. For $n=2$ and $D$ negative, we can solve the integral explicitly.
\[F(T)=\int_0^t\frac{dt}{1-Dy^2}=\tanh^{-1}(t\sqrt{D}).\]


We consider a similar question for pairs of primes. 

\begin{conjecture}
Let $p=f(a,b)$ and $q=f(a+c,b+d)$, where $p$ and $q$ are primes less than a large number $M$ and the slopes $b/a$ and $d/c$ are less than $T$ for  positive $a, b, c, d$.
If $f(a,b)=a^2+Db^2$  is a norm in the field $\Q(\sqrt{-D})$, then the distribution is given by
\[
D(\theta)
=\int_0^\theta\int_0^\theta\frac{(\alpha-\beta)-\sin(\alpha-\beta)\cos(\alpha-\beta)}{2\sin^3(\alpha-\beta)}d\alpha d\beta,
\]
where $\theta=\tan^{-1}(T\sqrt{D})$.
\end{conjecture}

For real quadratic fields the analogous conjecture is the following.
\begin{conjecture}
Let $p=f(a,b)$ and $q=f(a+c,b+d)$, where $p$ and $q$ are primes less than a large number $M$ and the slopes $b/a$ and $d/c$ are less than $T$ for  positive $a, b, c, d$.
If $f(a,b)=a^2-Db^2$ is a norm in the field $\Q(\sqrt{D})$  then the distribution is given by
\[
D(\theta)
=\int_0^\theta\int_0^\theta\frac{\sinh(\alpha-\beta)\cosh(\alpha-\beta)-(\alpha-\beta)}{2\sinh^3(\alpha-\beta)}d\alpha d\beta,
\]
where $\theta=\tanh^{-1}(T\sqrt{D})$.
\end{conjecture}

All of the conjectures can be expresses in in terms of multiple Dedekind zeta functions - both for primes and for pair of primes. In fact, this was the source of the conjectures. 
In Section 3, we introduce the version of multiple Dedekind zeta functions that we need. 
In Section 4, we use multiple Dedekind zeta functions, to find the right distributions. Also, in Section 4, we have given a summary of numerical tests of the distributions that we have obtained.

In the last Section, we compare sets primes in one number field to sets primes in another number field. That comparison is given in terms of remainder of a Dedekind zeta function at $s=1$. The comparison for distribution of pairs of primes is given by multiple residue of a multiple Dedekind zeta function at $s_1=1$ and $s_2=1$. Those comparisons are formulated as Conjectures, which are tested in many cases.

\section{Tests of Conjecture \ref{one prime} for primes in any number field}

\subsection{Cubic Fields}

Let $\alpha=\sqrt[n]{D}$. Let $\Z[\alpha]$ be the ring consisting of elements $a_0+a_1\alpha +\dots+a_{n-1}\alpha^{n-1}$, where $a_0, a_1,\dots,a_{n-1}$ are integers. We have that the norm $N(a+b\sqrt[n]{D})=a^n+(-1)^{n+1}Db^n$.

Let $n=3$ and $p=a^3+Db^3$ be a prime, for $D=\pm2,\pm3, \pm5, \pm 6$ Let $t=\frac{b}{a}$.  The distribution of primes with $p$ slope $t=\frac{b}{a}$ is $F^{-1}$, where
\[F(T)=\int_0^T\frac{dt}{
[
\sqrt[3]{1+Dt^3}
]^2
}.\]

For positive $D$, we tested the distribution of $b/a$ in the following way. 
We defined $t_1, t_2,\dots, t_9$ so that $F(t_1)=0.1,$ $F(t_2)=0.2, \dots, F(t_9)=0.9$. The approximation for $t_i$ that we used is up to the 4-th decimal point.

Another restriction that we impose is $b/a\leq 1$.
We might have that $t_9>1$ (or $t_8>1$) it that case the sequence $t_1,t_2,\dots$ might end at $t_8$ (or at $t_7$), respectively. 

For cubic fields we consider $N=4,000$ and $M=N^3/(1+D).$ Then, for primes less that $M$, we obtain

\[\frac{\max_k (C_{\sqrt[3]{-2}}(t_{k},t_{k+1}, M))}{\min_k(C_{\sqrt[3]{-2}}(t_{k},t_{k+1}, M))}=1.013.\]

\[\frac{\max_k (C_{\sqrt[3]{-3}}(t_{k},t_{k+1}, M))}{\min_k(C_{\sqrt[3]{-3}}(t_{k},t_{k+1}, M))}=1.015.\]

\[\frac{\max_k (C_{\sqrt[3]{-5}}(t_{k},t_{k+1}, M))}{\min_k(C_{\sqrt[3]{-5}}(t_{k},t_{k+1}, M))}=1.023.\]

\[\frac{\max_k (C_{\sqrt[3]{-6}}(t_{k},t_{k+1}, M))}{\min_k(C_{\sqrt[3]{-6}}(t_{k},t_{k+1}, M))}=1.022.\]



For negative $-D$, we consider primes of the form $a^3-Db^3$.
The values $t_k$ are defined as as an approximation of $F^{-1}(k/10)$ up to the 4th decimal, when that value is less that $0.1500$.
Then, for $N=4,000$ and $D=2,5,6$, and for primes less than for primes less that $M=N^n-N^n/(2^nD)$, we have that

\[\frac{\max_k (C_{\sqrt[3]{2}}(t_{k},t_{k+1}, M))}{\min_k(C_{\sqrt[3]{2}}(t_{k},t_{k+1}, M))}=1.011.\]
\[\frac{\max_k (C_{\sqrt[3]{5}}(t_{k},t_{k+1}, M))}{\min_k(C_{\sqrt[3]{5}}(t_{k},t_{k+1}, M))}=1.024.\]
\[\frac{\max_k (C_{\sqrt[3]{6}}(t_{k},t_{k+1}, M))}{\min_k(C_{\sqrt[3]{6}}(t_{k},t_{k+1}, M))}=1.032.\]

\subsection{A Quartic field}
Let $n=4$ and $p=a^4+11b^4$ be a prime. Let $t=\frac{b}{a}$.  The distribution for $t=\frac{b}{a}$ is $F^{-1}$, where 
\[F(t)=\int_0^t\frac{dy}{
[
\sqrt[4]{1+Dy^4}
]^2
}.\]

We define $t_k$ to be the approximation of $F^{-1}(k/10)$ up to the 4th decimal when $t_k<1$
For the quartic ring $\Z[\sqrt[4]{-11}]$, we take $N=3,500$ and $D=11$ and primes less than $M=N^4/(1+D)$. Then

\[\frac{\max_k (C_{\sqrt[4]{-11}}(t_{k},t_{k+1}, M))}{\min_k(C_{\sqrt[4]{-11}}(t_{k},t_{k+1}, M))}=1.0296.\]


\subsection{Degrees $5$ and $6$}

We include a few tests of the same distribution for polynomials of degree 5 and 6.

In the following examples, we consider intervals $(t_i,t_{i+1})$ for the slope $t=b/a$ such that $F(t_i)=\frac{i}{10}$ for $i=1,\dots,9$. Then we tested for primes less than a very large number $M$. For those primes we counted how many of the primes have slope $b/a$ in  each of the intervals $(t_i,t_{i+1})$ for $i=1,\dots,8$. We computed the standard deviation and the mean of those $8$ numbers. By error in these examples, we mean the standard deviation divided by the mean.

\[
\begin{tabular}{lll}
primes of the type &  less than $M$ & error=STDEV/mean\\
\hline\\
$p=a^5+2b^5$ & $1000^5 = 10^{15}$ & $0.028$\\
\\
$p=a^5+3a^4b+b^5$ & $750^5=2.3730469\times 10^{14}$ & $0.019$\\
\\
$p=a^6+3a^5b+b^6$ & $500^6=1.5625\times 10^{16}$ & $0.017$\\
\end{tabular}
\]

\subsection{Conjectures for distribution of primes in arbitrary number field}
We conjecture that if $\alpha$ is a solution of $f(x,1)=0$ where $f(x,y)=x^n+a_1x^{n-1}y+\dots+a_{n-1}xy^{n-1}+a_n$ is is an irreducible polynomial with integer coefficients, then for primes $p$ of the form $p=N(a+b\alpha)$  the statistics for $t=\frac{b}{a}$ has distribution
$t$ such that $F(t)$ is the constant distribution, where 
\[F(t)=\int_0^t\frac{dy}{[\sqrt[n]{f(1,y)}]^2}.\]

For $n=2$ and $D$ positive, we can solve the integral explicitly.
\[F(t)=\int_0^t\frac{dy}{1+Dy^2}=\tan^{-1}(t\sqrt{D}).\]
Then, we obtain exactly the distribution from Section 3.

For $n=2$ and $D$ negative, we can solve the integral explicitly.
\[F(t)=\int_0^t\frac{dy}{1-Dy^2}=\tanh^{-1}(t\sqrt{D}).\]

We state the following conjecture.

\begin{conjecture} 
\label{conj a prime}
Given a homogeneous polynomial $f(x,y)$ of degree $n$, we define
\[F(t)=\int_0^t\frac{dy}{[\sqrt[n]{f(1,y)}]^2},\]
and let $C_{f}(t_1,t_2,M)$ be the number of primes $p$ of the form $p=f(a,b)$ which are less that $M$, such that $F^{-1}(t_1)<b/a<F^{-1}(t_2)$. Consider the variables $t_1,t_2,t_3$ and $t_4$ in the interval $(0,T)$, where $F(t)$ is continuous. Let $t_1<t_2$ and $t_3<t_4$. Then,
\[\lim_{M\rightarrow \infty}\frac{C_{f}(t_1,t_2,M)}{C_{f}(t_3,t_4,M)}
=\frac{F(t_2)-F(t_1)}{F(t_4)-F(t_3)}.\]
\end{conjecture}

For quadratic fields this conjecture is a Theorem (see the Introduction together this section).
We have tested the conjecure for several cubic polynomials and a few of degrees $4$, $5$ and $6$; (earlier in this Section.)
Some of the data is available at \cite{Zouberou}.


\section{Multiple Dedekind Zeta Values and Related Functions}

In \cite{MDZF}, we defined multiple Dedekind zeta functions as certain integrals, they also have good expression in terms of infinite sum.
We will  recall the infinite sum that represents the multiple Dedekind zeta functions.

Let $K$ is a number field and let ${\cal{O}}_K$ be the ring of integers in the number field.  Let $C$ be a cone, which we will define soon.
\[\zeta_{K,C}(s_1,s_2)=\sum_{\gamma_1,\gamma_2 \in C}\frac{1}{N(\gamma_1)^{s_1}N(\gamma_1+\gamma_2)^{s_2}}.\]

More specifically, for $n=2$, we have that if $\gamma$ is in $C$ then $\gamma = a\mu_1 + b\mu_2$, where $a$ and $b$ are any positive integers and $\mu_1$ and $\mu_2$ are fixed totally positive algebraic integers in $K$, which means that $\sigma(\mu_i)>0$ for each real embedding $\sigma:K\rightarrow \R$ and $Re(\sigma(\mu_i))>0$ for each complex embedding $\sigma:K\rightarrow \C$.

So far this definition of multiple Dedekind zeta values coincides with the one in \cite{MDZF}.
We will make a few modifications of these definitions that will be useful for distribution of primes.

For a fixed real number $T>0$ we define $C_T$ to be
\[C_T=\{\gamma\in C \,\,\, | \,\,\,  \gamma = a\mu_1 + b\mu_2, \mbox{ where  } a \mbox{ and } b \mbox{ are {\it positive integers}, and } b/a<T \}.\]

We will define two new versions of multiple Dedekind zeta values: one discrete and one continuous. 

For the discrete one we use $C_T$ instead of $C$.
\[\zeta^{discrete}_{K,C_T}(s_1,s_2)=\sum_{\gamma_1,\gamma_2 \in C_T}\frac{1}{N(\gamma_1)^{s_1}N(\gamma_1+\gamma_2)^{s_2}}.\]

The other version is continuous. let 
\[f(a,b)=N(a\mu_1+b\mu_2)\]

If $\gamma_1=a\mu_1+b\mu_2$ and $\gamma_2 =c\mu_1+d\mu_2$, then
\[f(a+c,b+d) = N(\gamma_1 + \gamma_2).\]

With this notation we can express the discrete zeta as 
\[\zeta^{discrete}_{K,C_T}(s_1,s_2)=\sum_{a,b,c,d \geq 1;\,\, b/a<T; \,\,d/c<T}\frac{1}{f(a,b)^{s_1}f(a+c,b+d)^{s_2}}.\]

We can rewrite it in the form
\[\zeta^{discrete}_{K,C_T}(s_1,s_2)=\sum_{f(a,b) \geq 1;\,\, b/a<T; \,\,d/c<T}\frac{1}{f(a,b)^{s_1}f(a+c,b+d)^{s_2}}.\]

The continuous version is very similar. The key difference is that instead of summing we integrate with respect to a certain measure. We define 
\[\zeta^{continuous}_{K,C_T}(s_1,s_2)
=\int_{f(x,y)\geq 1;\,\,\, y/x<T; \,\,w/z<T}
\frac{1}{f(x,y)^{s_1}f(x+z,y+w)^{s_2}}
\omega(x,y)\omega(z,w)
\]
We are going to use also
\[\zeta^{continuous}_{K,C_T}(s)=\int_{f(x,y)\geq 1;\,\,\, y/x<T}
\frac{1}{f(x,y)^{s}}
\omega(x,y).\]

Where
$\omega(x,y)=\frac{dxdy}{f(x,y)^{1/n}}$ is a measure and $n$ is the degree of the homogeneous polynomial $f$. 


\section{Primes represented by polynomials and multiple Dedekind zeta functions}

\begin{conjecture}
Consider primes $p$ less than $M$ of the type $p=f(a,b)$, where $f$ is an irreducible homogeneous polynomial over $\Z$. 
Assume $a$ and $b$ are positive. For a given value of $T$ consider all $a$ and $b$ such that $f(a,b)$ is a prime and $0<b/a<T$. 
If we fix $M$ (large) and we let $T$ vary, then the distribution of the number of such primes is given by 
\[\zeta^{continuous}_{K,C_T}(2).\]
\end{conjecture}

\begin{lemma}
As a function of $T$, we have that
\[\zeta^{continuous}_{K,C_T}(2)=(constant)\cdot \int_0^T\frac{dt}{\left(\sqrt[n]{f(1,t)}\right)^2}.\]
\end{lemma}

\proof Consider the integral
\[\zeta^{continuous}_{K,C_T}(s)=\int_{f(x,y)\geq 1;\,\,\, y/x<T}\frac{1}{f(x,y)^{s}}\frac{dxdy}{f(x,y)}=\int_{f(x,y)\geq 1;\,\,\, y/x<t}\frac{1}{f(x,y)^{s+1}},\]
for positive $x$ and $y$.
We change the variables to $r=\sqrt[n]{(x,y)}$, where $n$ is the degree of the number field.
Let also $t=y/x$ be the slope with respect to which we consider the distribution. We want express $x$ and $y$ in terms of $r$ and $t$

We have $r^n=x^n f(1,y/x)=x^n f(1,t).$ Therefore 
\[x=\frac{r}{\sqrt[n]{f(1,t)}}.\]
Since $t=y/x$, we have $y=tx$ and 
\[y = \frac{rt}{\sqrt[n]{f(1,t)}}.\]
Now we compute the Jacobian for the change of variables.
\[\frac{\partial x}{\partial r} =  \frac{1}{\sqrt[n]{f(1,t)}}.\]
\[\frac{\partial y}{\partial r} =  \frac{t}{\sqrt[n]{f(1,t)}}.\]
\[\frac{\partial x}{\partial t} =  -\frac{1}{n}\cdot\frac{rf'_t(1,t)}{\sqrt[n]{f(1,t)}}.\]
\[\frac{\partial y}{\partial t} =  \frac{r}{\sqrt[n]{f(1,t)}}- \frac{1}{n}\cdot\frac{rtf'_t(1,t)}{f(1,t)^{\frac{n+1}{n}}}.\]
Therefore, $Jac=\frac{\partial x}{\partial r} \cdot \frac{\partial y}{\partial t} - \frac{\partial x}{\partial t} \cdot \frac{\partial y}{\partial r} = \frac{r}{\left(\sqrt[n]{f(1,t)}\right)^2}$.
We obtain that the zeta function
\begin{eqnarray}
&&\zeta^{continuous}_{K,C_T}(s)\\
&&=\int_{f(x,y)\geq 1;\,\,\, y/x<T}\frac{1}{f(x,y)^{s}}\omega(x,y)\\
&&=\int_{f(x,y)\geq 1;\,\,\, y/x<T}\frac{dxdy}{f(x,y)^{s+\frac{1}{n}}}\\
&&=\int_0^T\int_1^\infty\frac{1}{r^{ns+1}}\cdot \frac{drdt}{\left(\sqrt[n]{f(1,t)}\right)^2}\\
&&=\int_1^\infty\frac{dr}{r^{ns}}\int_0^T\frac{dt}{\left(\sqrt[n]{f(1,t)}\right)^2}
\end{eqnarray}
For any fixed value of $s$ the distribution 
is

\[F(T)=(constant)\cdot\int_0^T\frac{dt}{\left(\sqrt[n]{f(1,t)}\right)^2}.\]

\section{Pairs of primes and multiple Dedekind zeta functions}
\subsection{Gaussian integers}
\begin{conjecture}
Consider pairs primes $p$ and $q$ less that $M$ such that $p=N(\gamma_1)$ and $q=N(\gamma_1+\gamma_2)$ where $\gamma_1$ and $\gamma_2$ are Gaussian integers inside a sector with an angle $\theta$.
If we fix $M$ (large)  and we let $\theta$ vary,  the distribution is given by
\[D(\theta)=\int_0^\theta\int_0^\theta\frac{(\alpha-\beta)-\sin(\alpha-\beta)\cos(\alpha-\beta)}{2\sin^3(\alpha-\beta)}d\alpha d\beta.\]
which is a constant times $\zeta^{continuous}_{Q(i),C_T}(2,2)$, which depends on $T$, since $C_T$ does. The relation between $T$ and $\theta$ is $T=\tan(\theta)$
\end{conjecture}

{\bf Test of the distribution:} We did the following test for this distribution for Gaussian integers. Consider Gaussian integers $a+ib$ whose norm is a prime number less that $500^2$. For angles $\frac{k\pi}{36}$, we consider cones $C_k$ of Gaussian integers whose $\gamma$ such that $\arg(\gamma)$ is in the interval 
$\left(\frac{\pi}{36},\frac{(k+1)\pi}{36}\right)$. Consider any pair of Gaussian integers $\gamma_1$ and $\gamma_2$ that belong to the cone $C_k$, such that both $N(\gamma_1)$ and $N(\gamma_1+\gamma_2)$ are primes.
Let $P_k$ denote such pairs of primes for $k=1,2,3,4,5$.
We tested that the quotients $\frac{P_k}{P_l}$ are on ``average" one tenth of the percent of the prediction $D\left(\frac{k\pi}{36}\right)/D\left(\frac{l\pi}{36}\right)$. By ``average" we mean the following.
Let $q_k=\frac{P_k}{D\left(\frac{k\pi}{36}\right)}$. For $k=1,2,3,4,5$ compute the standard deviation (STDEV) divided by the mean of the sample. Think of the standard deviation as something measuring $maximum-minimum$. In fact, for a sample of $5$ elements, usually, the standard deviation  is three times less that the difference $maximum-minimum$. If we divide by the mean, we would measure by what fraction is the maximum bigger than the minimum times $1/3$.

We did a test  for primes pair of primes, each one less that $500^2$.
For the data $q_1,\dots,q_5$ we have $STDEV/mean=0.0008$, which is less than one tenth of a percent difference.

\begin{lemma}
We have that $\zeta^{continuous}_{Q(i),C_T}(2,2)$ is a constant times $D(\theta)$.
\end{lemma}
\proof
For the continuous multiple Dedekind zeta, we have
\begin{eqnarray*}
&&\zeta^{continuous}_{Q(i),C_T}(2,2) =\\
&&=\int_{\Delta}\frac{1}{f(x,y)^{s_1}f(x+z,y+w)^{s_2}}\omega(x,y)\omega(z,w)=\\
&&= \int_\Delta\frac{1}{(x^2+y^2)^{2}((x+z)^2+(y+w)^2)^{2}}\frac{dxdy}{\sqrt{x^2+y^2}}\frac{dzdw}{\sqrt{z^2+w^2}},
\end{eqnarray*}
where $f(x,y)=x^2+y^2$,
$\omega(x,y)=\frac{dxdy}{\sqrt{x^2+y^2}}$, 
and
\[\Delta=\{(x,y,z,w)\,|;\,\,x>0;\,\,y>0,\,\,z>0; \,\,w>0;\,\,0<y/x<T;\,\, 0<w/z<T;\,\,x^2+y^2>1\}\]
We change the Cartesian coordinates $(x,y)$ to polar $(r,\alpha)$ and similarly the Cartesian coordinates $(z,w)$ to polar $(R,\beta)$ such that $x+iy=re^{i\alpha}$ and $z+iw=Re^{i\beta}$. Then 
\[\omega(x,y)=\frac{dxdy}{\sqrt{x^2+y^2}}=drd\alpha\] and 
\[\omega(z,w)=\frac{dzdw}{\sqrt{z^2+w^2}}=dRd\beta.\]

Then \[(x+z)^2+(y+w)^2=(r\cos\l\alpha+R\cos\beta)^2+(r\sin\alpha+R\sin\beta)^2=r^2+R^2+2rR\cos(\alpha-\beta).\]
Let the slope $T=\tan\theta$.
We obtain that the continuous zeta at $(2,2)$ is 
\[\zeta^{continuous}_{Q(i),C_T}(2,2) =\int_{r>1;\,\,R>0;\,\, 0<\alpha<\theta;\,\, 0<\beta<\theta}\frac{1}{(r^2)^2(r^2+R^2+2rR\cos(\alpha-\beta))^2}dr d\alpha dR d\beta.\]

Let $s=R/r$. Then $drdR=rdrds$. Therefore,
\begin{eqnarray}
&&\int_{r>1;\,\,s>0;\,\, 0<\alpha<\theta;\,\,0<\beta<\theta}\frac{1}{(r^2)^2\left(r^2+r^2s^2+2r^2s\cos(\alpha-\beta)\right)^2} rdr ds d\alpha d\beta\\
&&= \int_1^\infty \frac{dr}{r^7}\int_0^\theta\int_0^\theta\left(\int_0^\infty\frac{ds}{(s^2+2s\cos(\alpha-\beta)+1)^2}\right)d\alpha d\beta.
\end{eqnarray}

We are going to use that 
\[
\int \frac{dv}{(v^2+1)^2}=\frac{1}{2}\tan^{-1}v+\frac{1}{2}\frac{v}{v^2+1}.\]
Let $u=s+\cos(\alpha-\beta)$ and let
$v=u/\sin(\alpha-\beta)$. Assume that $\alpha>\beta$. Then
\begin{eqnarray}
&&\int_0^\infty\frac{ds}{(s^2+2s\cos(\alpha-\beta)+1)^2}\\
&&=
\int_{\cos(\alpha-\beta)}\frac{du}{(u^2+\sin^2(\alpha-\beta))^2}\\
&&=
\int_{\cot(\alpha-\beta)}^\infty\frac{\sin(\alpha-\beta)dv}{(v^2\sin^2(\alpha-\beta)+\sin^2(\alpha-\beta))^2}\\
&&=
\frac{1}{\sin^3(\alpha-\beta)}\int_{\cot(\alpha-\beta)}^\infty\frac{dv}{(v^2+1)^2}\\
&&=
\frac{1}{\sin^3(\alpha-\beta)}\left[\frac{1}{2}\tan^{-1}v+\frac{1}{2}\frac{v}{v^2+1}\right]_{\cot(\alpha-\beta)}^\infty\\
&&=
\frac{1}{\sin^3(\alpha-\beta)}\left[\frac{\pi}{4}-\frac{1}{2}\left(\frac{\pi}{2}-(\alpha-\beta)\right)-\frac{1}{2}\sin(\alpha-\beta)\cos(\alpha-\beta)\right]\\
&&=
\frac{(\alpha-\beta)-\sin(\alpha-\beta)\cos(\alpha-\beta)}{2\sin^3(\alpha-\beta)}.
\end{eqnarray}
Therefore,
\begin{eqnarray}
&&\zeta^{continuous}_{Q(i),C_t}(2,2)=\\
&&= \int_1^\infty \frac{dr}{r^7}\int_0^\theta\int_0^\theta\left(\int_0^\infty\frac{ds}{(s^2+2s\cos(\alpha-\beta)+1)^2}\right)d\alpha d\beta=\\
&&=\int_1^\infty\frac{dr}{r^7}\int_0^\theta\int_0^\theta\frac{(\alpha-\beta)-\sin(\alpha-\beta)\cos(\alpha-\beta)}{2\sin^3(\alpha-\beta)}d\alpha d\beta=\\
&&=\frac{1}{8}D(\theta).
\end{eqnarray}

We used $D(\theta)$ to test the conjecture, which makes it computationally more effective compared to
the integral defining $\zeta^{continuous}_{\Q(\sqrt{-d}),C_T}(2,2)$.

\subsection{Imaginary Quadratic Fields}

The formulation of the corresponding conjecture is almost identical to the one for Gaussian integers. 
\begin{conjecture}
Let $K=\Q(\sqrt{-d})$ be an imaginary quadratic field and let $N:K\rightarrow \Q$ be the norm map.
Consider pairs primes $p$ and $q$ less that $M$ such that $p=N(\gamma_1)$ and $q=N(\gamma_1+\gamma_2)$ where $\gamma_1$ and $\gamma_2$ are algebraic integers in the field $K$ inside a sector with an angle $\theta$.
If we fix $M$ (large)  and we let $\theta$ vary,  the distribution is given by
\[D(\theta)=\int_0^\theta\int_0^\theta\frac{(\alpha-\beta)-\sin(\alpha-\beta)\cos(\alpha-\beta)}{2\sin^3(\alpha-\beta)}d\alpha d\beta.\]
which is a constant times $\zeta^{continuous}_{K,C_T}(2,2)$, where $T=\tan(\theta)$
\end{conjecture}

\begin{lemma}
The function $\zeta^{continuous}_{K,C_T}(2,2)$ is the same form every imaginary quadratic field $K$.
\end{lemma}
\proof
Let $x+y\sqrt{-d}=re^{i\alpha}$ and $z+w\sqrt{-d}=Re^{i\beta}$. Then, the slopes are $\tan(\alpha)=y\sqrt{d}/x$ and $\tan(\beta)=w\sqrt{d}/z$
If $y'=y\sqrt(d)$ and $w'=w\sqrt{d}$, then $\tan(\alpha)=y'/x$, $\tan(\beta)=w'/z$, $x^2+dy^2 = x^2+(y')^2$ and $z^2+dw^2 = z^2+(w')^2$
Therefore,
\begin{eqnarray}
&&\zeta^{continuous}_{\Q(\sqrt{-d}),C_T}(2,2)\\
&&=\zeta^{continuous}_{\Q(i),C_T}(2,2)\\
&&= (constant)\cdot\int_0^\theta\int_0^\theta\frac{(\alpha-\beta)-\sin(\alpha-\beta)\cos(\alpha-\beta)}{2\sin^3(\alpha-\beta)}d\alpha d\beta\\
&&=(constant)\cdot D(\theta)
\end{eqnarray}

{\bf Test of the distributions:}
For the field $K=Q(\sqrt{-d})$ with $d<10$, we computed the number $P_k(-d)$ of pairs primes $(p,q)$ less than $M$ (different choice of the constant $M$ for different field) such that $p=N(\gamma_1)$ and $q=N(\gamma_1+\gamma_2)$, where
$\gamma_1$ and $\gamma_2$ 
are algebraic integers in $K$ and both $\arg(\gamma_1)$ and $\arg(\gamma_2)$ are in the interval $\left(\frac{\pi}{36}, \frac{(k+1)\pi}{36}\right)$. 
The following is a table of quadratic extensions $Q(\sqrt{-d})$ and the corresponding upper bounds $M$ for the primes $p$ and $q$.
\[
\begin{tabular}{lll}
$Q(\sqrt{-d})$ & pairs of primes less than & error of prediction\\
\hline\\
$Q(\sqrt{-1})$ & $250,000$ &  $0.0008$\\
\\
$Q(\sqrt{-2})$ & $500,000$ & $0.006$\\
\\
$Q(\sqrt{-3})$ & $750,000$ & $0.008$\\
\\
$Q(\sqrt{-5})$ & $1,250,000$ & $0.006$\\
\\
$Q(\sqrt{-6})$ & $1,500,000$ & $0.095$\\
\\
$Q(\sqrt{-7})$ & $1,750,000$ & $0.003$\\
\end{tabular}
\]
Then again
the proportions $P_1(-d):P_2(-d):\dots:P_5(-d)$ are almost the same as $D\left(\frac{\pi}{36}\right):D\left(\frac{2\pi}{36}\right):\dots:D\left(\frac{5\pi}{36}\right)$. Again, the error is very small.
More precisely, if $q_k(-d)=\frac{P_k(-d)}{D\left(\frac{k\pi}{36}\right)}$, then for the sample $q_1(-d),\dots,q_5(-d)$, we have that $STDEV/mean$ is less than $1\%$ for each number field $\Q(\sqrt{-d})$ with $d<10$.
 
 \subsection{Real Quadratic Fields}

\begin{conjecture}
Let $K$ be a real quadratic field and let $N:K\rightarrow \Q$ be the norm map.
Consider pairs primes $p$ and $q$ less that $M$ such that $p=N(\gamma_1)$ and $q=N(\gamma_1+\gamma_2)$ where $\gamma_1$ and $\gamma_2$ are algebraic integers of $K$ inside a sector with a hyperbolic angle 
$\theta$.
If we fix $M$ (large)  and we let $\theta$ vary,  the distribution is given by
which is a constant times $\zeta^{continuous}_{K,C_t}(2,2)$, where $t=\tanh(\theta)$
\end{conjecture}
{\bf Remark:} By a hyperbolic angle we mean the following. The hyperbolic angle between the vectors $(r\cosh(\alpha),r\sinh(\alpha))$ and $(R\cosh(\beta),R\sinh(\beta))$ is $\alpha-\beta$, where $r$ and $R$ are any positive constants. 

\begin{lemma}
The function $\zeta^{continuous}_{K,C_T}(2,2)$ is the same form every real quadratic field $K$.
It is equal to a constant multiple of
\[D_{real}(\theta)=\int_0^\theta\int_0^\theta\frac{\sinh(\alpha-\beta)\cosh(\alpha-\beta)-(\alpha-\beta)}{2\sinh^3(\alpha-\beta)}d\alpha d\beta\]

\end{lemma}
\proof
Let $(x,y\sqrt{d})=(r\cosh(\alpha),r\sinh(\alpha))$ and $z+w\sqrt{d}=(R\cosh(\beta),R\sinh(\beta))$. Then the hyperbolic  slopes are $\tanh(\alpha)=y\sqrt{d}/x$ and $\tanh(\beta)=w\sqrt{d}/z$
If $y'=y\sqrt(d)$ and $w'=w\sqrt{d}$, then $\tanh(\alpha)=y'/x$, $\tanh(\beta)=w'/z$, $x^2-dy^2 = x^2-(y')^2$ and $z^2-dw^2 = z^2-(w')^2$
Similarly, we obtain
\begin{eqnarray}
&&\zeta^{continuous}_{\Q(\sqrt{d}),C_t}(2,2)\\
&&= (constant)\cdot\int_0^\theta\int_0^\theta\frac{\sinh(\alpha-\beta)\cosh(\alpha-\beta)-(\alpha-\beta)}{2\sinh^3(\alpha-\beta)}d\alpha d\beta\\
&&=(constant)\cdot D_{real}(\theta)
\end{eqnarray}
Instead of computing the integral from the very beginning, we can use that pass from trigonometric functions to hyperbolic ones by $\sin(i\theta)=i\sinh(\theta)$ and $\cos(i\theta)=\cosh(\theta)$.

{\bf Test of the distributions:}
For several field $K=Q(\sqrt{d})$ with $d<14$, we computed the number $P_k(d)$ of pairs primes $(p,q)$ less than $M$ such that $p=N(\gamma_1)$ and $q=N(\gamma_1+\gamma_2)$, where
$\gamma_1$ and $\gamma_2$ 
are algebraic integers in $K$ and both $\tanh^{-1}(\gamma_1)$ and $\tanh^{-1}(\gamma_2)$ are in the interval $(\pi/36, (k+1)\pi/36)$. Then again
 the proportions $P_1(d):P_2(d):\dots:P_5(d)$ are almost the same as $D(\pi/36):D(2\pi/26):\dots:D(5\pi/36)$. Again, the error is very small.
 More precisely, if $q_k(d)=P_k(d)/D(k\pi/36)$, then for the sample $q_1(d),\dots,q_5(d)$, we have that $STDEV/mean$ is less than one percent for each number field $\Q(\sqrt{d})$.
 
 \[
\begin{tabular}{lll}
$Q(\sqrt{d})$ & pairs of primes less than & error of prediction\\
\hline\\
$Q(\sqrt{2})$ & $100,000$ & $0.004$\\
\\
$Q(\sqrt{3})$ & $100,000$ & $0.008$\\
\\
$Q(\sqrt{5})$ & $100,000$ & $0.004$\\
\\
$Q(\sqrt{13})$ & $156,250$ & $0.008$\\
\end{tabular}
\]

\section{Comparing the number of primes from pairs of number fields}

From the previous section it is clear that the distribution up to a constant that described totally split primes in an imaginary quadratic field is independent of the choice of the imaginary quadratic field. 
Similarly, for primes in a real quadratic field - the distribution is independent of the choice of a real quadratic field. 

For pair of primes in an imaginary quadratic field we have one distribution up to a constant. And pairs of primes in a real quadratic field we have a different distribution up to a constant.

In this section, we state several conjectures on comparison between sets of primes in two imaginary quadratic fields or two real quadratic fields.

Let $C_1$ be the cone whose elements have angle $\theta$ between $\pi/36$ and $2\pi/36$. That is, if $t=b/a$ then $\theta=\tan^{-1}t$ in the imaginary quadratic case, and $\theta=\tanh^{-1}t$ in the real quadratic case.

\begin{conjecture}
Let $K=Q(\sqrt{-d})$ be an imaginary number field. Consider primes less that $M$ that split completely whose factors lie in the cone $C_1$. Denote this number by $P_M(-d)$. Let $h(-d)$ is the class number for that quadratic field.
Then 
\[\lim_{M\rightarrow\infty}\frac{P_M(-d_1)}{P_M(-d_2)} = \frac{h(-d_2)}{h(-d_1)}.\]
\end{conjecture}

We have shown this experimentally for all $d<20$ and primes less that $1,000^2$ we have errors less than one percent.
 
\begin{conjecture}
Let $K=Q(\sqrt{d})$ be a real quadratic field. Consider primes less that $M$ that split completely whose factors lie in the cone $C_1$. Denote this number by $P_M(d)$. Let $h(d)$ be the class number and
and $R(d)$ the  regulator of the number field. Then 
\[\lim_{M\rightarrow\infty}\frac{P_M(d_1)}{P_M(d_2)} = \frac{h(d_2)R(d_2)}{h(d_1)R(d_1)}.\] 
\end{conjecture}

{\bf Test of the conjecture:}
\[
\begin{tabular}{lll}
\hline\\

discriminant & (number of primes in $C_1$)$\times$(regulator of the number field)\\
\hline\\
$d:76$ & $23947.9225771526$\\
$d:60$ & $12012.017518793384$\\
$d:17$ & $24054.53632226616$\\
$d:56$ & $24058.379117099787$\\
$d:44$ & $24027.144008192117$\\
$d:40$ & $12012.657309687034$\\
$d:28$ & $24003.773460713535$\\
$d:24$ & $24042.18973884143$\\
$d:12$ & $24058.546031358073$\\
$d:13$ & $24078.17173282505$\\
$d:8$ & $24084.895388125315$\\
\hline\\
\end{tabular}
\]

We have shown this experimentally for all $d<20$ and primes less that $1,000^2$ we have errors less than 1.2 percent.

\begin{conjecture}
\label{conj P-d}
Let $Q(\sqrt{-d})$ be an imaginary number field. Consider a pairs  of primes $(p,q)$ less that $M$ that split completely. Assume $p=N(\gamma_1)$  and $q=N(\gamma_1+\gamma_2)$, where $\gamma_1$ and $\gamma_2$ lie in the cone $C_1$. Denote this number by $PQ_N(-d)$. Let $h(-d)$ be the class number for that quadratic field.
\[\lim_{M\rightarrow\infty}\frac{PQ_M(-d_1)}{PQ_M(-d_2)} = \left(\frac{h(-d_2)}{h(-d_1)}\right)^2.\]

\end{conjecture}
 We have shown this experimentally for all $d<14$ and primes less that $500^2$ we have errors less than 2 percent.

\begin{conjecture}
\label{conj Pd}
Let $K=Q(\sqrt{d})$ be an imaginary number field. Consider a pairs  of primes $(p,q)$ less that $M$ that split completely. Assume $p=N(\gamma_1)$  and $q=N(\gamma_1+\gamma_2)$, where $\gamma_1$ and $\gamma_2$ lie in the cone $C_1$. Denote this number by $PQ_N(d)$. Let $h(d)$ be the class number and $R(d)$ be the regulator of the number field. Then
\[\lim_{M\rightarrow\infty}\frac{PQ_M(d_1)}{PQ_M(d_2)} = \left(\frac{h(d_2)R(d_2)}{h(d_1)R(d_1)}\right)^2.\] 

\end{conjecture}

\[
\begin{tabular}{lll}
discriminant & primes less than  $25,000,000$ in $C_1$\\
\hline\\
$d:8$ & $P_{25,000,000}(8)=27202$\\
$d:12$ & $P_{25,000,000}(12)=18139$\\
\end{tabular}
\]

\[
\begin{tabular}{lll}
discriminant  & pair of primes, each less than $1,000$ in $C_1$\\
\hline\\
$d:8$ & $PQ_{1,000}(8)=39829$\\
$d:12$ & $PQ_{1,000}(12)=17650$\\
\end{tabular}
\]
\[\left(\frac{P_{25,000,000}(8)}{P_{25,000,000}(12)}\right)^2\div\left(\frac{PQ_{1,000}(8)}{PQ_{1,000}(12)}\right)=\frac{39829 \times18139^2}{17650\times 27202^2}=1.0034\]

{\bf Remark:}
1. One can show that for a real quadratic fields the multiple Dedekind zeta values times a suitable power of the square root of the discriminant is a period of an unramified mixed Tate motive over the corresponding number ring \cite{periods}.
 
2. One can think of the quotients $\frac{h(-d_1)}{h(-d_2)}$ or $\frac{h(d_1)R(d_1)}{h(d_2)R(d_2)}$ as \[\frac{\sqrt{d_1}Res_{s=1}\zeta_{\Q(\sqrt{d_1})}(s)}{\sqrt{d_2}Res_{s=1}\zeta_{\Q(\sqrt{d_2})}(s)},\] where $d_1$ and $d_2$ are discriminant of quadratic number fields (both negative or both positive)

A reformulation of Conjecture \ref{conj P-d}, which is test with a good precision, is the following.
\begin{conjecture}
For imaginary quadratic fields, we consider primes less that $M$ that split completely whose factors lie in the cone $C_1$. Denote this number by $P_N(-d)$. Let $h(-d)$ is the class number for that quadratic field.
Then 
\[\lim_{M\rightarrow\infty}\frac{P_M(-d_1)}{P_M(-d_2)} = \frac{\sqrt{d_1}Res_{s=1}\zeta_{\Q(\sqrt{-d_1})}(s)}{\sqrt{d_2}Res_{s=1}\zeta_{\Q(\sqrt{-d_2})}(s)}.\]
\end{conjecture}

A reformulation of Conjecture \ref{conj Pd}, which is test with a good precision, is the following.
\begin{conjecture}
For real quadratic fields, we consider primes less that $M$ that split completely whose factors lie in the cone $C_1$. Denote this number by $P_N(-d)$. Let $h(-d)$ is the class number for that quadratic field.
Then 
\[\lim_{M\rightarrow\infty}\frac{P_M(d_1)}{P_M(d_2)} = \frac{\sqrt{d_1}Res_{s=1}\zeta_{\Q(\sqrt{d_1})}(s)}{\sqrt{d_2}Res_{s=1}\zeta_{\Q(\sqrt{d_2})}(s)}.\]
\end{conjecture}

{\bf Remark:} 
1. For pairs of primes in a quadratic field, the conjectures are related it to a multiple Dedekind zeta function. Moreover, when comparing different fields the quotient is given by a residue of a multiple Dedekind zeta function. 
Note that $\left(\frac{h(-d_1)}{h(-d_2)}\right)^2$ or $\left(\frac{h(d_1)R(d_1)}{h(d_2)R(d_2)}\right)^2$ is equal to 
 \[\left(\frac{(\sqrt{d_1})Res_{s=1}\zeta_{\Q(\sqrt{d_1})}(s)}{(\sqrt{d_2})^2Res_{s=1}\zeta_{\Q(\sqrt{d_2})}(s)}\right)^2,\] when both $d_1$ and $d_2$ are negative or both a positive.

2. We recall a formula for a double residues of a certain multiple Dedekind zeta function \cite{MDZF}. For the double residue of that function we have
\[Res_{s_1=1}Res_{s_2=1}\zeta_{K,C}(s_1,s_2) = [Res_{s=1}\zeta_{K,C}(s)]^2,\]
where 
\[\zeta_{K,C}(s)=\sum_{\gamma \in C}\frac{1}{N(\gamma)^{s}}.\]

We would expect
\[\lim_{M\rightarrow\infty}\frac{PQ_M(-d_1)}{PQ_M(-d_2)} = \frac{(\sqrt{d_1})^2Res_{s_1=1}Res_{s_2=1}\zeta_{\Q(\sqrt{-d_1})}(s_1,s_2)}{(\sqrt{d_2})^2Res_{s_1=1}Res_{s_2=1}\zeta_{\Q(\sqrt{-d_2})}(s_1,s_2)}.\]

Let $C$ be a cone for a real quadratic number $K$, which means the following. For fixed totally positive 
algebraic integers $\mu$ and $\nu$ in $K$.  We consider all elements $\gamma= a\mu +b\nu$, where $a$ and $b$ are positive integers. Then,
\[C=\{\gamma \,\,   |  \,\, \gamma = a\mu+ b\nu\}.\] Then by definition

\[\zeta_{K,C}(s_1,s_2)=\sum_{\gamma_1,\gamma_2 \in C}\frac{1}{N(\gamma_1)^{s_1}N(\gamma_1+\gamma_2)^{s_2}}.\]

We would expect that for a suitable definition of multiple Dedekind zeta values $\zeta_K(s_1,s_2)$ we would like to have
\begin{equation}
\label{residue}
Res_{s_1=1}Res_{s_2=1}\zeta_{K}(s_1,s_2) = [Res_{s=1}\zeta_{K}(s)]^2,
\end{equation}
where
$Res_{s=1}\zeta_{K}(s)$ is the classical Dedekind zeta function. The for the field $\Q(\sqrt{3})$ we have that  Equation \ref{residue} is true. 
Take for instance a cone $C$ consisting of all numbers $a +b (2+\sqrt{3})$, where $a$ is a positive integer and $b$ is zero or a positive integer, we have
that the multiple Dedekind zeta function $\zeta_{\Q(\sqrt{3}),C}(s)$ coincides with  the Dedekind zeta function $\zeta_{\Q(\sqrt{3}),C}(s)$. We can take
$\zeta_{\Q(\sqrt{3})}(s_1,s_2)$ to be $\zeta_{\Q(\sqrt{3}),C}(s_1,s_2)$ for that particular cone $C$. With that definition of $\zeta_{\Q(\sqrt{3})}(s_1,s_2)$, we have that Equation \ref{residue} holds.

For other real quadratic number field, there is a suggestion by Zagier.  One might need to take a finite union of small cones, in order to make the multiple Dedekind zeta function independent of the cones. Then, follow the techniques in the paper \cite{MDZF} in order to compute the double residues.

\section*{Acknowledgements:} Ivan Horozov would like to thank the Max Planck Institute in Bonn for the great working conditions and the financial support. He is also very thankful for a very useful and detailed discussion with Don Zagier.

Nickola Horozov would like to thank his mentor Louis Cooper for introducing him to Python.

Zouberou Sayibou would like to thank the LSAMP program at City University of New York - BCC  for the opportunity and financial support to do research.

\renewcommand{\em}{\textrm}

\begin{small}
   
\end{small}

\section*{Affiliations:}

\hspace{.45cm} Ivan Horozov, City University of New York - BCC and the Graduate Center;  

ivan.horozov@bcc.cuny.edu

Nickola Horozov, Farragut Middle School;
alokcin.sci@gmail.com

Zouberou Sayibou, 
Stanford University; zouberou@stanford.edu

\end{document}